\documentclass[11pt]{article}  
\usepackage{amsfonts,amsmath,amssymb,latexsym}
\usepackage{bbm}
\usepackage[dvips]{graphicx,epsfig}
\usepackage[latin1]{inputenc}

\def\real{\mathbb{R}}

\def\H{\mathbb{H}}

\def\mat1{\mathbb{I}}

\def\dzero{\mbox{{\scriptsize $\mathbb{O}$}}}  
\def\dun{\mathbbm 1}

\newcommand { \Iletter}[1] {I\kern-0.10em #1 }

\def\bit{\begin{itemize}}
\def\eit{\end{itemize}}
\def\ben{\begin{enumerate}}
\def\een{\end{enumerate}}
\def\bde{\begin{description}}
\def\ede{\end{description}}

\def\bar{\begin{array}}
\def\ear{\end{array}}
\def\beq{\begin{equation}}
\def\eeq{\end{equation}}
\def\bfi{\begin{figure}[hbt] \begin{center}}
\def\efi{\end{center} \end{figure}}
\def\noi{\noindent}
\def\bce{\begin{center}}
\def\ece{\end{center}}
\newcommand{\proof}{{\bf Proof. }}
\newcommand{\cqfd}{\hfill $\Box$}

\newtheorem {theo} {Theorem}[section]

\newtheorem {coro} {Corollary} [section]

\newtheorem {lemm} {Lemma}[section]

\newtheorem {propo} {Proposition}[section]

\newtheorem {rema} {Remark} [section]

\begin{document}

\title{The Fr\'echet Contingency Array Problem is Max-Plus Linear}
\author{L. TRUFFET \\
Ecole des Mines de Nantes\\
4, rue A. Kastler, La Chantrerie \\
BP 20722, Nantes 44307 Cedex 3 France.\\
e-mail: Laurent.Truffet@emn.fr
}

\maketitle 

\begin{abstract}
In this paper we show that the so-called array Fr\'echet problem 
in Probability/Statistics is $(\max, +)$-linear. 
The upper bound of Fr\'echet is obtained using simple 
arguments from residuation theory and lattice distributivity. The 
lower bound is obtained as a loop invariant of a greedy 
algorithm. The algorithm is based on the max-plus linearity 
of the Fr\'echet problem and the Monge property of 
bivariate distribution.
\end{abstract}

\noi

\section{Introduction}
As a preliminary remark the author would like this paper 
to be a modest tribute to 
the work of Maurice Fr\'echet in Statistics. The work 
has been started at the occasion of his $130$th 
birthday and the $100$th anniversary of his stay in Nantes 
(the town the author is living in) as professor in 
Mathematics.

In this paper it is shown that the tropical or 
max-plus semiring $\real_{\max}$ (i.e. 
the set $\real$ of real numbers with  $\max$ as addition 
and $+$ as multiplication, see the precise definition in 
subsection~\ref{subRmax}) 
is the underlying algebraic structure which 
is well suited to the Fr\'echet contingency (or 
correlation) array problem \cite{kn:Frechet51}. 
In other words the Fr\'echet problem is a tropical problem 
which thus has its place in the new trends of 
idempotent mathematics 
founded by V. P. Maslov and its collaborators in the 
1980s (see e.g.\cite{kn:Litvinov06} and references therein).

From this main result the Fr\'echet upper bound is 
derived by residuation and the distributivity property 
of $\real_{\max}$ as a lattice. The Fr\'echet lower 
bound is obtained as a loop invariant of a greedy 
algorithm. This algorithm is based on the tropical 
nature of the Fr\'echet problem and the Monge property 
of bivariate distribution.

\subsection{Organization of the paper}

The paper is written to be sef-contained. Thus, in 
Section~\ref{secPrel} we introduce main notations 
used in the paper, we restate the Fr\'echet array 
problem and its bounds. We define the tropical 
semiring $\real_{\max}$ and recall basic 
results on residuation theory and lattices. In 
Section~\ref{secMR} we prove the main result of the 
paper that is the Fr\'echet array problem is 
max-plus linear in the space of cumulative distribution 
functions (see Theorem~\ref{thMR}). From this result in 
Section~\ref{secNewApproach} we derive the Fr\'echet bounds 
using new approaches. The upper bound is derived from 
residuation theory and the lattice distributivity 
property of the max-plus semiring $\real_{\max}$ 
(see Corollary~\ref{corUp}). The lower bound 
is obtained as the loop invariant of a 
greedy algorithm (see Proposition~\ref{propLow}). 
In Section~\ref{secConcl} we conclude this work.

\section{Preliminaries}
\label{secPrel}

In this Section we recall basic results concerning 
Fr\'echet array problem and the max-plus semi-ring 
$\real_{\max}$.

\subsection{The Fr\'echet contingency array problem and its solution}

This problem is described in e.g. \cite{kn:Frechet51}. Let $n$ be an
integer $\geq 1$. The set $\textsf{Mat}_{nm}(\real_{+})$ denotes the set of $n \times m$
matrices which entries are nonnegative real numbers. We define the
partial order $\overset{\cal{D}}{\preceq}$ on
$\textsf{Mat}_{nn}(\real_{+})$ as follows:
\begin{equation}
\label{compDistrib2D}
A=[a_{i,j}] \overset{\cal{D}}{\preceq} B=[b_{i,j}] \overset{\textsf{def}}{\Leftrightarrow} 
\forall i,j, \; 
\sum_{l=1}^{i} \sum_{k=1}^{j} a_{l,k} \leq \sum_{l=1}^{i} \sum_{k=1}^{j} b_{l,k}.
\end{equation}
Introducing the fundamental $n \times n$ matrix:
\begin{equation}
\label{matD}
D \stackrel{\textsf{def}}{=} [1_{\{i \leq j \}}],
\end{equation}
where $1_{\{i \leq j \}} =1$ if $i \leq j$ and $0$ otherwise, the 
partial order $\overset{\cal{D}}{\preceq}$ can be rewritten as 
follows:
\begin{equation}
A \overset{\cal{D}}{\preceq} B \Leftrightarrow D A D^{T} \leq D B D^{T} \mbox{ (entrywise)}
\end{equation}  
where $()^{T}$ denotes the transpose operator.

Let $p,q \in \textsf{Mat}_{n1}(\real_{+})$ such that 
$\sum_{i=1}^{n} p_{i}= \sigma = \sum_{j=1}^{n} q_{j}$. Without 
loss of generality we can assume:
\[
\sigma =1.
\]
The problem of Fr\'echet is then to find (if exist) the 
maximum and the minimum element w.r.t the partial 
order $\overset{\cal{D}}{\preceq}$ of the subset of 
$\textsf{Mat}_{nn}(\real_{+})$: 
\begin{equation}
\label{ensFrechet}
\mathcal{H}(p,q) \stackrel{\textsf{def}}{=} \{ F \in \textsf{Mat}_{nn}(\real_{+}) | 
F \mbox{ satisfies ({\bf F1})  and ({\bf F2})} \}
\end{equation}
with: \\
({\bf F1}). $F \underline{1} = p$, \\
({\bf F2}). $ \underline{1}^{T} F = q^{T}$. \\
Where $\underline{1}$ denotes the $n$-dimensional vector which coordinates 
are $1$'s. 

Fr\'echet proved that there exist a maximum element, 
$\bigvee_{\overset{\cal{D}}{\preceq}} \mathcal{H}(p,q) \stackrel{\textsf{not.}}{=} 
F_{\max}$, and a 
minimum element, $\bigwedge_{\overset{\cal{D}}{\preceq}} 
\mathcal{H}(p,q) \stackrel{\textsf{not.}}{=} F_{\min}$ such that:
\begin{subequations}
\begin{equation}
\label{defFrechMax}
(D F_{\max} D^{T})_{i,j} = \min((D p)_{i}, (q^{T} D^{T})_{j}) 
\stackrel{\textsf{not.}}{=} (\overline{F}_{\max})_{i,j},
\end{equation}
and 
\begin{equation}
\label{defFrechMin}
(D F_{\min} D^{T})_{i,j} = \max(0, (D p)_{i} + (q^{T} D^{T})_{j} - \sigma) 
\stackrel{\textsf{not.}}{=} (\overline{F}_{\min})_{i,j}
\end{equation}
\end{subequations}
for all $i,j=1, \ldots, n$.

\subsection{The max-plus semiring $\real_{\max}$}
\label{subRmax}

Let $\real$ be the field of real numbers. The max-plus semiring 
$\real_{\max}$ is the set $\real \cup \{-\infty\}$ equipped with 
the addition $\oplus : \; (a,b) \mapsto a \oplus b 
\stackrel{\textsf{def}}{=} \max(a,b)$ 
and the multiplication $\odot : \; (a,b) \mapsto a \odot b 
\stackrel{\textsf{def}}{=} a + b$. 
The neutral element for $\oplus$ is $\dzero := -\infty$ and the 
neutral element for $\odot$ is $\dun := 0$. 
The max-plus semiring is said to be an idempotent semiring because 
the addition is idempotent, i.e. $a \oplus a =a$, $\forall a$.   
 
The idempotent semiring $\real_{\max}$ or idempotent semirings
isomorphic to $\real_{\max}$ has many applications in discrete
mathematics, algebraic geometry, computer science, computer languages,
linguistic problems, optimization theory, discrete event systems,
fuzzy theory (see e.g. \cite{kn:Bac-cooq},
\cite{kn:Cuninghame-Green79}, \cite{kn:Golan99}, \cite{kn:RGetal03}).

\subsection{Order properties of $\real_{\max}$ and residuation}

Let us consider the max-plus semiring $\real_{\max}$ already 
defined in the introduction. The binary relation $\mathcal{R}$ 
defined by: $a \mathcal{R} b \overset{\textsf{def}}{\Leftrightarrow} a \oplus b 
=b$ coincides with the standard order $\leq$ on $\real$. We denote 
$\overline{\real}_{\max}$ 
the semiring completed by adjoining to $\real_{\max}$ a $\top:=+\infty$ 
element which satisfies $a \oplus \top = \top$, $\forall a$, 
$a \odot \top = \top$, $\forall a \neq \dzero$, and $\dzero \odot \top = \dzero$. 
This is mainly motivated by the fact that some of the further results can be stated 
in a simpler way in $\overline{\real}_{\max}$.

The completed max-plus semiring $\overline{\real}_{\max}$ 
is a complete sup-semilattice, 
i.e. $\forall A \subseteq \overline{\real}$ $\oplus A \stackrel{\textsf{def}}{=} 
\underset{x \in A}{\oplus} x$ exists 
in $\overline{\real} \stackrel{\textsf{def}}{=}\real \cup\{-\infty, +\infty\}$. This 
implies that $\overline{\real}_{\max}$ 
is also a complete 
inf-semilattice because $\wedge A = \oplus \{x \in \overline{\real} | 
\forall a \in A, \; x \leq a \}$. Thus, $\overline{\real}_{\max}$ is a 
complete lattice. Finally, let us mention that $\overline{\real}_{\max}$ is 
a distributive lattice, i.e.: 
\begin{equation}
\label{lattDistrib}
a \oplus (b \wedge c) = (a \wedge b) \oplus (a \wedge c).
\end{equation}
As we will see in the sequel this property will be of particular 
importance. Let us define left and right division 
in $\overline{\real}_{\max}$ by: 
$b / a \stackrel{\textsf{def}} = 
\oplus \{x \in \overline{\real} | x \odot a \leq b\}$ and 
$a \backslash b \stackrel{\textsf{def}}{=} \oplus \{x \in \overline{\real} | 
a \odot x  \leq b\}$. Left and right division coincide with the 
usual subtraction to which we add the following properties:
$\dzero \backslash a = \top$, $\top \backslash a = \dzero$ if $a \neq \top$, $\top$ 
otherwise (similar formulae for $/$).

We extend operations and binary relations from scalars to matrices as follows. 
If $A=[a_{i,j}]$, $B=[b_{i,j}]$ then: $A \leq B \mbox{ (entrywise) } \Leftrightarrow 
\forall i,j, \; a_{i,j} \leq b_{i,j}$, $A \oplus B = [a_{i,j} \oplus b_{i,j}]$, 
$A \wedge B = [a_{i,j} \wedge b_{i,j}]$, and $A \odot B$ denotes the 
matrix such that its entry $(i,j)$, $(A \odot B)_{i,j}$, is: 
$(A \odot B)_{i,j} = \oplus_{k} a_{i,k} \odot b_{k,j}$. We obviously have: 
$(A \odot B)^{T} = B^{T} \odot A^{T}$. We also extend the divisions to 
(possibly rectangular) matrices with suitable dimension:
\begin{subequations}
\begin{equation}
\label{AsousB}
(A \backslash B)_{i,j}  \stackrel{\textsf{def}}{=} (\oplus \{X | A \odot X \leq B\})_{i,j} = 
\wedge_{k} a_{k,i} \ b_{k,j}, 
\end{equation}
\begin{equation}
\label{DsurC}
(D / C)_{i,j} \stackrel{\textsf{def}}{=}(\oplus\{X | X \odot C \leq D \})_{i,j} = 
\wedge_{l} d_{i,l} / c_{j,l}. 
\end{equation}
\end{subequations}
The $\oplus$ in the formulae (\ref{AsousB}) and (\ref{DsurC})
corresponds to the supremum w.r.t entrywise order between
matrices. The application $Y \mapsto A \backslash Y$ (resp. $Y \mapsto
Y / C$) is called the {\em residuated mapping} of the application $X
\mapsto A \odot X$ (resp. $X \mapsto X \odot C$). 

For more details on max-plus algebra and residuation 
theory we refer the reader to e.g \cite[Chp. 4 and references 
therein]{kn:Bac-cooq}.

\section{Main Result}
\label{secMR}

We begin this Section by the following fundamental lemma.
\begin{lemm}
\label{lem1}
Let $(u^{j})_{j=1}^{m}$ be $m$ elements of $\textsf{Mat}_{n1}(\real_{+})$. 
Define $s^{j} = \sum_{k=1}^{j} u^{k}$. Then, 
\begin{equation}
\label{eqlem1}
[ u^{1} \cdots u^{n}] \; \underline{1} = [s^{1} \cdots s^{n}] \odot \underline{\dun} 
\end{equation}
where $\underline{\dun}$ denotes the $n$-dimensional vector which all components 
are $\dun$'s.
\end{lemm}
\proof Let us remark that: $[ u^{1} \cdots u^{n}] \; \underline{1} = s^{n}$. 
Because the vectors $u^{j}$ have all their coordinates nonnegative: 
$s^{1} \leq \cdots \leq s^{n}$ (componentwise), which is equivalent 
to: $s^{n} = s^{1} \oplus \cdots \oplus s^{n}$. Now, we just 
have to remark that: $s^{1} \oplus \cdots \oplus s^{n}=
[s^{1} \cdots s^{n}] \odot \underline{\dun}$ which ends the proof 
of the result. \cqfd

\begin{rema}
\label{remsurlem1}
In the previous Lemma define $U=[u^{1} \cdots u^{n}]$ and $S = [s^{1} \cdots 
s^{n}]$. We remark that $S = U D^{T}$, recalling that $D$ is the matrix 
defined by (\ref{matD}). Thus, relation (\ref{eqlem1}) can be rewritten: 
\begin{equation}
U \underline{1} = (U D^{T}) \odot \underline{\dun}.
\end{equation}
\end{rema}

\begin{theo}[Main Result]
\label{thMR}
Let us consider a matrix $F \in \textsf{Mat}_{nn}(\real_{+})$. 
Let $p$ and $q$ be two elements of $\textsf{Mat}_{n1}(\real_{+})$. 
Then, 
\[
F \in \mathcal{H}(p,q)
\Longleftrightarrow
\left\{ \bar{ll} (D F D^{T}) \odot \underline{\dun} & = D p \\
                 \underline{\dun}^{T} \odot (D F D^{T})   & = q^{T} D^{T}
\ear \right.
\]
\end{theo}
\proof $F \in \mathcal{H}(p,q) \Leftrightarrow \left\{ \bar{ll} F \underline{1} & = p \\
                 \underline{1}^{T} F & = q^{T} 
\ear \right.$. \\
Because matrix $D$ is invertible so is $D^{T}$ and: 
\[
\left\{ \bar{ll} F \underline{1} & = p \\
                 \underline{1}^{T} F & = q^{T} 
\ear \right.
\Longleftrightarrow
\left\{ \bar{lll} D F \underline{1} & = D p & \mbox{ (eq 1)}\\
                 \underline{1}^{T} F D^{T}   & = q^{T} D^{T} & \mbox{ (eq 2)}.
\ear \right.
\]
For (eq 1). We apply Lemma~\ref{lem1} with $u^{j}:= (D F)_{., j}$ and 
$s^{j} = (D F D^{T})_{.,j}$ be the $j$th column vectors of 
matrices $D F$ and $DFD^{T}$, respectively. We obtain the following 
equality: $ D p = D F \underline{1} =(D F D^{T}) \odot \underline{\dun}$. \\
For (eq 2). We apply Lemma~\ref{lem1} 
with $u^{j} := (D F^{T})_{.,j}$ and $s^{j} := (D F^{T} D^{T})_{.,j}$. 
We have: $D F^{T} \underline{1}= (q^{T} D^{T})^{T} = (D F^{T} D^{T}) 
\odot \underline{\dun}$. 
By definition of $()^{T}$ we have: $q^{T} D^{T} = 
\underline{\dun}^{T} \odot 
(D F D^{T})$ which ends the proof of the Theorem. \cqfd

This result can be reformulated as follows. Let us define the
following sets: 
\begin{equation}
\label{defDistn1}
\textsf{Dist}_{n1} \stackrel{\textsf{def}}{=} D
\;\textsf{Mat}_{n1}(\real_{+}) = \{ D x, \; x \in \textsf{Mat}_{n1}(\real_{+}) \}
\end{equation}
and 
\begin{equation}
\label{defDistnn}
\textsf{Dist}_{nn}
\stackrel{\textsf{def}}{=} D \;\textsf{Mat}_{nn}(\real_{+}) \; D^{T} = 
\{ D X D^{T}, \; X \in \textsf{Mat}_{nn}(\real_{+}) \}
\end{equation}
and 
for all $\overline{P}, \overline{Q} \in \textsf{Dist}_{n1}$:
\begin{equation}
\H(\overline{P}, \overline{Q}) \stackrel{\textsf{def}}{=} \{ \overline{F} \in 
\textsf{Dist}_{nn} | \overline{F} \odot \underline{\dun} = \overline{P} \mbox{ and } 
\underline{\dun}^{T} \odot \overline{F} = \overline{Q}^{T} \}.
\end{equation}
Then, Theorem~\ref{thMR} states that $\forall p, q \in \textsf{Mat}_{n1}(\real_{+})$ 
and $\forall F \in \textsf{Mat}_{nn}(\real_{+})$:
\begin{equation}
\label{equivMR}
F \in \mathcal{H}(p,q) \Leftrightarrow D F D^{T} \in \H(D p, D q).
\end{equation}
Or, equivalently:
\begin{equation}
\label{equivMR2}
D \; \mathcal{H}(p,q) \; D^{T} = \H(D p, D q).
\end{equation}

\section{New approach for the Fr\'echet bounds}
\label{secNewApproach}

From our main result, Theorem~\ref{thMR}, we obtain Fr\'echet bounds 
by methods which seem to be new to the best knowledge of 
the author.

\subsection{Upper bound}

The Fr\'echet upper bound is obtained as a direct consequence of Theorem~\ref{thMR}.

\begin{coro}[Fr\'echet upper bound]
\label{corUp}
Let  $p$ and $q$ be two elements of $\textsf{Mat}_{n1}(\real_{+})$.
Under the condition that $p^{T} \underline{1} = q^{T} \underline{1}$ the 
set $\mathcal{H}(p,q)$ is not empty and the upper Fr\'echet bound 
$F_{\max}$ is such that $D F_{\max} D^{T} \stackrel{\textsf{not.}}{=} \overline{F}_{\max}$ 
is the greatest sub-solution of the following max-plus linear system of equations: 
\begin{equation}
\label{syscoro}
\left\{
\bar{ll} \overline{F} \odot \underline{\dun} & = D p \\
\underline{\dun}^{T} \odot \overline{F} & = q^{T} D^{T}
\ear \right.
\end{equation} 
that is:
\[
\overline{F}_{\max}= 
((D p) / \underline{\dun}) \wedge (\underline{\dun}^{T} \backslash (q^{T} D^{T})).
\]
\end{coro}
\proof Let us study (\ref{syscoro}) when replacing $=$ by $\leq$. Then, 
by definition of $/$ we have: $\overline{F} \odot \underline{\dun} \leq D p 
\Leftrightarrow \overline{F} \leq (D p) / \underline{\dun}$. And by definition 
of $\backslash$ we have: $\overline{F} \leq \underline{\dun}^{T} \backslash 
(q^{T} D^{T})$. The two previous inequalities are equivalent 
to: $\overline{F} \leq ((D p) / \underline{\dun}) \wedge (\underline{\dun}^{T} \backslash 
(q^{T} D^{T})) = \overline{F}_{\max}$. Now, we have to prove that 
(A): $\overline{F}_{\max} \odot \underline{\dun} = D p$ and (B): 
$\underline{\dun}^{T} \odot \overline{F}_{\max} = q^{T} D^{T}$. \\
Let us prove (A). \\
For all $i =1, \ldots, n$ we write:
\[
\bar{ll}
(\overline{F}_{\max} \odot \underline{\dun})_{i} & = \oplus_{j=1}^{n} ((D p) / 
\underline{\dun}) \wedge (\underline{\dun}^{T} \backslash (q^{T} D^{T}))_{i,j} \odot \dun \\
\mbox{} & = \oplus_{j=1}^{n} (Dp)_{i} \wedge (q^{T} D^{T})_{j} \\
\mbox{} & = (D p)_{i} \wedge (\oplus_{j=1}^{n} (q^{T} D^{T})_{j}) 
\mbox{  \bigg(by lattice distributivity (\ref{lattDistrib}) \bigg)} \\
\mbox{} & = (D p)_{i} \wedge (q^{T} D^{T})_{n} \mbox{  \bigg($\forall j$: 
$(q^{T} D^{T})_{j} 
\leq (q^{T} D^{T})_{n}$ \bigg)} \\
\mbox{} & = (D p)_{i} \wedge (D p)_{n} \mbox{  \bigg($(b^{T} D^{T})_{n} = q^{T} 
\underline{1} = p^{T} \underline{1} = (D p)_{n}$ \bigg)}  \\
\mbox{} & = (D p)_{i} \mbox{  \bigg($\forall i$: $(D p)_{i} \leq (D p)_{n}$ \bigg)}.
\ear
\]
We prove (B) similarly. Hence the result is now achieved. \cqfd

\subsection{Lower bound}
\label{subsecLow}
In this section we obtain Fr\'echet lower bound by a greedy 
algorithm based on max-plus linearity of the 
Fr\'echet problem and the well-known Monge property (see e.g. \cite{kn:Burkardetal96}) of elements 
of the set $\textsf{Dist}_{nn}$ defined by (\ref{defDistnn}), that 
is for all $\overline{F} \in \textsf{Dist}_{nn}$: \\
\noi 
({\bf M}). $\forall i,j=0, \ldots n-1$: $\overline{F}_{i,j} \odot
\overline{F}_{i+1,j+1} \geq \overline{F}_{i,j+1} \odot \overline{F}_{i+1,j}$, \\
with the convention that $\forall k$ $\overline{F}_{0,k} = \overline{F}_{0,k} = 0$. \\
\noi
Let $p, q \in  \textsf{Mat}_{n1}(\real_{+})$ be two vectors such that 
$p^{T} \; \underline{1} 
= q^{T} \; \underline{1}=1$. Let us consider the following algorithm.

\noi
\texttt{Lower}(n,p,q) \\
$\forall i =1, \ldots , n$, $\overline{F}_{i,n} := (D p)_{i}$ ; (a) \\
$\forall j =1, \ldots , n$, $\overline{F}_{n,j} := (q^{T} D^{T})_{j}$ ; (b) \\
\noi
For $j=n-1$ to $1$ do \\

For $i=n-1$ to $1$ do \\
\[
\overline{F}_{i,j} := \overline{F}_{i+1,j+1}^{-1} \odot (\overline{F}_{i,j+1} \odot 
\overline{F}_{i+1,j}) \oplus \dun
\]

end \\

\noi
end.

\begin{propo}
\label{propLow}
The algorithm $\texttt{Lower}$ computes the Lower bound of the 
Fr\'echet contingency array problem.
\end{propo}
\proof The initial conditions (a) and (b) of the algorithm $\texttt{Lower}$ 
come from the max-plus linearity of the Fr\'echet problem and that the 
Monge property ({\bf M}) implies:
\[
\forall i \leq i', j \leq j', \; \overline{F}_{i,j} \leq \overline{F}_{i',j'}.
\]
The proof is obtained by recurrence (see e.g. the detailed proof of this result 
by Fr\'echet himself \cite[p. 13]{kn:Frechet60}) \\

\noi
Denoting $\alpha_{l}= (D p)_{l}$, $\beta_{k} = (q^{T} D^{T})_{k}$ we have to 
prove that the loop invariant of the algorithm \texttt{Lower} corresponds 
to the Fr\'echet lower bound, i.e. $\forall l,k$: $\overline{F}_{l,k} = 1^{-1} 
\odot \alpha_{l} \odot \beta_{k} \oplus \dun$. \\

It is easy to see that the previous relation is true for $l=n$ with $k=1, \ldots , n$ 
and for $l=1, \ldots, n$ with $k=n$. Now, let us assume that the loop invariant 
is true for $(k,l) \geq (i,j)$, $(k,l) \neq (i,j)$. We have:
\[
\bar{ll}
\overline{F}_{i,j+1} \odot \overline{F}_{i+1,j} & = (1^{-1} \odot \alpha_{i} \odot \beta_{j+1} 
\oplus \dun ) \odot (1^{-1} \odot \alpha_{i+1} \odot \beta_{j} \oplus \dun) \\
\mbox{} &= 2^{-1} \odot \alpha_{i} \odot \alpha_{i+1} \odot  \beta_{j} \odot \beta_{j+1} \\
\mbox{} &\oplus 1^{-1} \alpha_{i} \odot \beta_{j+1} \oplus 1^{-1} \odot \alpha_{i+1} 
\odot \beta_{j} \oplus \dun.
\ear
\]
Because $(\real_{+}, \leq)$ is a totally ordered lattice: 
\[
\overline{F}_{i+1,j+1}= 1^{-1} \odot \alpha_{i+1} \odot \beta_{j+1} \oplus 
\dun \; \in \{1^{-1} \odot \alpha_{i+1} \odot \beta_{j+1}, \dun\}
\]
Thus, we have two cases to study. \\

\noi
\underline{$1$rst case}: $\overline{F}_{i+1,j+1} = 1^{-1} \odot \alpha_{i+1} \odot \beta_{j+1}$ \\
Let us compute:
\[
\bar{ll}
\overline{F}_{i+1,j+1}^{-1} \odot \overline{F}_{i,j+1} \odot \overline{F}_{i+1,j} & = 1^{-1} \odot \alpha_{i}
 \odot \beta_{j} \oplus \alpha_{i} \odot \alpha_{i+1}^{-1} \oplus \beta_{j} \odot \beta_{j+1}^{-1} \\
\mbox {} & \oplus 1 \odot  (\alpha_{i+1} \odot \beta_{j+1})^{-1}.
\ear
\]
Then, we just have to remark that: 
$\alpha_{i} \odot \alpha_{i+1}^{-1} = p_{i+1}^{-1} \leq \dun$, 
$\beta_{j} \odot \beta_{j+1}^{-1} = q_{j+1}^{-1} \leq \dun$ and 
$\overline{F}_{i+1,j+1} = 1^{-1} \odot \alpha_{i+1} \odot \beta_{j+1} \Leftrightarrow 
1 \odot  (\alpha_{i+1} \odot \beta_{j+1})^{-1} \leq \dun$. Thus, 
\[
\overline{F}_{i,j} = \overline{F}_{i+1,j+1}^{-1} \odot \overline{F}_{i,j+1} \odot \overline{F}_{i+1,j} 
\oplus \dun = 1^{-1} \odot \alpha_{i} \odot \beta_{j} \oplus \dun.
\] 

\noi
\underline{$2$nd case}: $\overline{F}_{i+1,j+1} = \dun$ \\
\[
\bar{ll}
\overline{F}_{i+1,j+1}^{-1} \odot \overline{F}_{i,j+1} \odot \overline{F}_{i+1,j} & = 
\overline{F}_{i,j+1} \odot \overline{F}_{i+1,j} \\
\mbox{} & = 
2^{-1} \odot \alpha_{i} \odot \alpha_{i+1} \odot  \beta_{j} \odot \beta_{j+1}
\oplus 1^{-1} \alpha_{i} \odot \beta_{j+1} \\
\mbox{} & \oplus 1^{-1} \odot \alpha_{i+1} \odot \beta_{j}
\oplus \dun 
\ear
\]
which could be rewritten as follows:
\[
\overline{F}_{i,j} = 1^{-1} \odot \alpha_{i} \odot \beta_{j} 
\odot (1^{-1} \odot \alpha_{i+1} \odot \beta_{j+1} 
\oplus q_{j+1} \oplus p_{i+1}) \oplus \dun.
\]
We remark that $q_{j+1}, p_{i+1} \geq \dun$. 
We also note that $\overline{F}_{i+1,j+1} = \dun \Rightarrow 
1^{-1} \odot \alpha_{i+1} \odot \beta_{j+1} \leq \dun$. Thus, 
we deduce because $\odot$ is non-decreasing and the definition of 
the standard order $\leq$ that:
\[
\overline{F}_{i,j} \geq 1^{-1} \odot \alpha_{i} \odot \beta_{j} \oplus \dun.
\] 

On the other hand $\overline{F}_{i+1,j+1} = \dun \Rightarrow 1^{-1} \odot \alpha_{i} 
\odot \beta_{j} \leq p_{i+1}^{-1} \odot q_{j+1}^{-1}$. And we deduce that:
\[
\bar{ll}
\overline{F}_{i,j} & \leq  p_{i+1}^{-1} \odot q_{j+1}^{-1} (1^{-1} \odot 
\alpha_{i+1} \odot \beta_{j+1} 
\oplus q_{j+1} \oplus p_{i+1} ) \oplus \dun \\
\mbox{}  & = 1^{-1} \odot \alpha_{i} \odot \beta_{j} \oplus p_{i+1}^{-1} \oplus q_{j+1}^{-1} \oplus \dun \\
\mbox{}  & = 1^{-1} \odot \alpha_{i} \odot \beta_{j} \oplus 
\dun \mbox{ (because $p_{i+1}^{-1}, q_{j+1}^{-1} \leq \dun$)}
\ear
\]
We conclude because $\leq$ is antisymmetric.

\cqfd

\section{Conclusion}
\label{secConcl}

In this paper we proved that the Fr\'echet correlation array problem 
is max-plus linear in the space of cumulative distribution function 
$\textsf{Dist}_{nn}$ defined by (\ref{defDistnn}). This remark leads to 
new methods to obtain the Fr\'echet bounds.

As a further work we would like to extend results of the 
paper to the continuous case based on the remark that : 
\[
\int_{\real} f(x,y) dy = \sup_{z \in \real} \bigg(\int_{-\infty}^{z} f(x,y) dy \bigg)
\]
for all nonnegative functions $f$ such that $\int_{\real} f(x,y) dy$ exists.

\end{document}